\documentclass[a4paper,11pt,reqno]{amsart}
\usepackage{amsmath, amsthm, amssymb, mathrsfs, amsfonts}
\usepackage{array}
\usepackage[active]{srcltx}
\usepackage{graphicx}
\usepackage[margin=3cm]{geometry}
\usepackage{cite}
\usepackage[colorlinks=true,allcolors=blue]{hyperref}
\usepackage{enumerate}
\usepackage[square,numbers]{natbib}

\newtheorem{teor}{Theorem}[section]

\newtheorem{lem}{Lemma}[section]

\newtheorem{cor}{Corollary}[section]

\newtheorem{obs}{Remark}
\newtheorem{conj}{Conjecture}[section]

\def\car{\mathop{\rm char}\nolimits}

\usepackage{tikz,xcolor,hyperref}

\definecolor{lime}{HTML}{A6CE39}
\DeclareRobustCommand{\orcidicon}{%
	\begin{tikzpicture}
	\draw[lime, fill=lime] (0,0)
	circle [radius=0.16]
	node[white] {{\fontfamily{qag}\selectfont \tiny ID}};
	\draw[white, fill=white] (-0.0625,0.095)
	circle [radius=0.007];
	\end{tikzpicture}
	\hspace{-2mm}
}

\foreach \x in {A, ..., Z}{%
	\expandafter\xdef\csname orcid\x\endcsname{\noexpand\href{https://orcid.org/\csname orcidauthor\x\endcsname}{\noexpand\orcidicon}}
}

\title[Group identities on symmetric units \& oriented involutions]%
      {Group identities on symmetric units under \\ \ oriented involutions in group algebras}

\author[A. Holgu\'in-Villa]{Alexander Holgu\'in-Villa\orcidA{}}
\address{Alexander Holgu\'in-Villa, Escuela de Matem\'aticas, Universidad Industrial de Santander}
\email{aholguin@uis.edu.co}

\author[J. H. Castillo]{John H. Castillo\orcidB{}}
\address{John H. Castillo, Departamento de Matem\'aticas y Estad\'istica, Universidad de Nari\~no}
\email{jhcastillo@udenar.edu.co}

\begin{document}
\maketitle
 \noindent
  \renewcommand{\refname}{References}

%%%%%%%%%%%%%%%%%%%%%%%%%%%%%%%%%%%%%%%%%%%%%%%%%%%%%%%%%%%%%%%%%%

\begin{abstract}
  Let $\mathbb{F}G$ denote the group algebra of a locally finite group $G$ over the infinite field $\mathbb{F}$
  with $\car(\mathbb{F})\neq 2$, and let $\circledast:\mathbb{F}G\rightarrow \mathbb{F}G$ denote the involution
  defined by $\alpha=\Sigma\alpha_{g}g \mapsto \alpha^\circledast=\Sigma\alpha_{g}\sigma(g)g^{\ast}$, where
  $\sigma:G\rightarrow \{\pm1\}$ is a group homomorphism (called an orientation) and $\ast$ is an involution of
  the group $G$. In this paper we prove, under some assumptions, that if  the $\circledast$-symmetric units of
  $\mathbb{F}G$ satisfies a group identity then $\mathbb{F}G$ satisfies a polynomial identity, i.e., we give an
  affirmative answer to a Conjecture of B. Hartley in this setting. Moreover, in the case when the prime radical
  $\eta(\mathbb{F}G)$ of $\mathbb{F}G$ is nilpotent we characterize the groups for which the symmetric units
  $\mathcal{U}^+(\mathbb{F}G)$ do satisfy a group identity.
\end{abstract}

%%%%%%%%%%%%%%%%%%%%%%%%%%%%%%%%%%%%%%%%%%%%%%%%%%%%%%%%%%%%%%%%%%%%%%%%%%
\section{Introduction}

Let $\mathbb{F}G$ denote the group algebra of the group $G$ over the field $\mathbb{F}$. Any involution
$\ast: G \rightarrow G$ can be extended $\mathbb{F}$-linearly to an algebra involution of $\mathbb{F}G$.
Such a map is called a group involution of $\mathbb{F}G$. A natural example is the so-called \emph{classical
involution}, which is induced from the map $g\mapsto g^\ast=g^{-1}$, for all $g\in G$.

Let $\sigma:G\rightarrow \{\pm1\}$ be a non-trivial homomorphism (called an \emph{orientation} of $G$).
If $\ast:G \rightarrow G$ is a group involution, an \emph{oriented group involution} of $\mathbb{F}G$ is
defined by

\begin{equation}\label{eq0}
\alpha=\sum_{g\in G} \alpha_gg \mapsto \alpha^\circledast=\sum_{g\in G} \alpha_g\sigma(g)g^{\ast}.
\end{equation}

Notice that, as $\sigma$ is non-trivial, $\car(\mathbb{F})$ must be different from $2$. It is clear that,
$\alpha\mapsto \alpha^\circledast$ is an involution in $\mathbb{F}G$ if and only if $gg^\ast\in N=ker(\sigma)$
for all $g\in G$.

In the case when the involution on $G$ is the classical involution, the map $\circledast$ is precisely the oriented
involution introduced by S.~P. Novikov, \cite{Nov:70}, in the context of $K$-theory.

We denote with $\mathbb{F}G^+=\{\alpha\in \mathbb{F}G:\alpha^\circledast=\alpha\}$  and
$\mathbb{F}G^-=\{\alpha\in \mathbb{F}G: \alpha^\circledast=-\alpha\}$ the sets of symmetric and skew-symmetric elements of
$\mathbb{F}G$ under $\circledast$ and, writing $\mathcal{U}(\mathbb{F}G)$ for the group of units of $\mathbb{F}G$,
we let $\mathcal{U}^+(\mathbb{F}G)$ denote the set of symmetric units, i.e.,
$\mathcal{U}^+(\mathbb{F}G)=\{\alpha\in \mathcal{U}(\mathbb{F}G): \alpha^\circledast=\alpha\}$.

Let $\langle x_1, x_2,\ldots\rangle$ be the free group on a countable set of generators. If $H$ is any subset
of a group $G$, we say that $H$ satisfies a group identity ($H\in$ GI or $H$ is GI for short) if there exists
a non-trivial reduced word $\omega(x_1, x_2, \ldots, x_n)\in \langle x_1, x_2, \ldots\rangle$ such that
$\omega(h_1, h_2, \ldots, h_n)=1$ for all $h_i\in H$. For instance, if we write $(x_1, x_2)=x_1^{-1}x_2^{-1}x_1x_2$
and $(x_1, x_2, \ldots, x_n, x_{n+1})=((x_1, x_2, \ldots, x_n), x_{n+1})$, for all $n\geq 2$, then $\langle H\rangle$
is abelian if it satisfies the group identity $(x_1, x_2)$, nilpotent if it satisfies $(x_1, x_2, \ldots, x_n)$,
for some $n$ and $n$-Engel if it satisfies $(x_1, \underbrace{x_2, x_2, \ldots, x_2}_{n \text{ times}})$ for some $n$.

Some time ago and with the idea of establishing a connection between the additive and multiplicative structure of a group
algebra $\mathbb{F}G$, Brian Hartley made the following famous conjecture:

\begin{conj}[Hartley's Conjecture]
Let $G$ be a torsion group and $\mathbb{F}$ a field. If the unit group $\mathcal{U}(\mathbb{F}G)$ of $\mathbb{F}G$
satisfies a group identity, then $\mathbb{F}G$ satisfies a polynomial identity.
\end{conj}

%A conjecture of Brian Hartley states that if the unit group $\mathcal{U}(\mathbb{F}G)$ of the group algebra
%$\mathbb{F}G$ of a torsion group $G$ over a field $\mathbb{F}$ satisfies a group identity, then $\mathbb{F}G$
%satisfies a polynomial identity.

Let $R$ be an $\mathbb{F}$-algebra. Recall that a subset $S$ of $R$ satisfies a polynomial identity ($S\in \text{PI}$ or
$S$ is PI for short) if there exists a non-zero polynomial $f(x_1, x_2, \ldots, x_n)$ in the free associative algebra
$\mathbb{F}\{X\}$ on the countable infinite set of non-commuting variables $X=\{x_1, x_2, \ldots\}$ such that
$f(a_1, \ldots, a_n)=0$ for all $a_i\in S$. For instance, $R$ is commutative if it satisfies the polynomial identity
$f(x_1, x_2)=x_1x_2-x_2x_1$ and, any finite dimensional associative algebra satisfies the \emph{standard polynomial identity}
of degree $n+1$, where $n=\dim_{\mathbb{F}}R$ \cite[Lemma 5.1.6, p. 173]{Pas:77},
$$
St_{n+1}(x_1, x_2, \ldots, x_{n+1})=\sum_{\rho\in \mathcal{S}_{n+1}}(sgn\rho)x_{\rho(1)}x_{\rho(2)}\cdot\cdot\cdot x_{\rho(n+1)}.
$$
Group algebras $\mathbb{F}G$ satisfying a PI were classified in two subsequent papers of Passman and
Isaacs-Passman, see \cite[Corollaries 5.3.8 and 5.3.10, p. 196-197]{Pas:77}.

Giambruno, Jespers and Valenti \cite{GJV:94} solved the Hartley's conjecture for semiprime group rings, and Giambruno,
Sehgal and Valenti \cite{GSV:97} solved it in general for group algebras over infinite fields. By using the results
of \cite{GSV:97}, Passman \cite{Pas:97} gave necessary and sufficient conditions for $\mathcal{U}(\mathbb{F}G)$ to
satisfy a group identity, when $\mathbb{F}$ is infinite. Subsequently, Liu \cite{Liu:99} confirmed that the conjecture
also holds for finite fields and Liu and Passman in \cite{LP:99} extended the results of \cite{Pas:97} to this case.
The same question for groups with elements of infinite order was studied by Giambruno, Sehgal and Valenti in \cite{GSV:00}.
For further details about these results see Lee \cite[Chapter 1]{Lee:10} and the references quoted therein.

Let $\ast$ be an involution of a group algebra $\mathbb{F}G$ induced by an involution of the group $G$, the so-called group involution.
When $\ast$ is the classical involution induced from $g\mapsto g^{-1}, g\in G$, Giambruno, Sehgal and Valenti \cite{GSV:98}
showed that if $G$ is a torsion group, $\mathbb{F}$ is infinite with $\car(\mathbb{F})\neq 2$, and $\mathcal{U}^+(\mathbb{F}G)$
satisfies a group identity then $\mathbb{F}G$ satisfies a polynomial identity. They also classified groups $G$ such that
$\mathcal{U}^+(\mathbb{F}G)$ satisfies a group identity. Sehgal and Valenti \cite{SV:06} extended the results of \cite{GSV:98}
to non-torsion groups under the usual restriction for the only if part related to Kaplansky's Conjecture
(the units of $\mathbb{F}G$ are trivial if $G$ is a torsion-free group and $\mathbb{F}$ is a field).

Considering group involutions $\ast$, i.e., $\ast$ is an involution on $G$ extended $\mathbb{F}$-linearly to the group algebra
$\mathbb{F}G$, Dooms and Ruiz \cite{DRM:07} proved the following.

\begin{teor}[{\cite[Theorem 3.1]{DRM:07}}]\label{teor01}
Let $\mathbb{F}$ be an infinite field with $\car(\mathbb{F})\neq 2$ and let $G$ be a non-abelian group such that
$\mathbb{F}G$ is regular. Let $\ast$ be an involution on $G$. Suppose one of the following conditions holds:

\begin{itemize}
 \item[\emph{(i)}] $\mathbb{F}$ is uncountable,
 \item[\emph{(ii)}] all finite non-abelian subgroups of $G$ which are $\ast$-invariant have no simple components
                   in their group algebra over $\mathbb{F}$ that are non-commutative division algebras other than
                   quaternion algebras.
\end{itemize}
Then $\mathcal{U}^+(\mathbb{F}G)\in$ GI if and only if $G$ is an SLC-group with canonical
involution given by the expression \eqref{eq1} below. Moreover, in this case $\mathbb{F}G^+$ is a ring contained in $\zeta(\mathbb{F}G)$.
\end{teor}

Using the last result and under some assumptions, Dooms and Ruiz proved that if $\mathcal{U}^+(\mathbb{F}G)$ is
GI then $\mathbb{F}G$ is PI, giving an affirmative answer to the Hartley's Conjecture in this setting.
They also characterized, with mild restrictions, the locally finite groups for which the symmetric units
$\mathcal{U}^+(\mathbb{F}G)$ satisfy a group identity, when the prime radical $\eta(\mathbb{F}G)$ of $\mathbb{F}G$ is nilpotent.
Giambruno, Polcino Milies and Sehgal \cite{GPS:09i} completely solved the question for  group algebras of torsion groups, with
group involutions such that $\mathcal{U}^+(\mathbb{F}G)$ is GI.

In the classification results on group algebras whose symmetric units with respect to the classical involution satisfy a group
identity in some sense the exceptional cases turned out to involve Hamiltonian $2$-groups, \cite[Theorem 7, p. 459]{GSV:98},
because they are non-abelian groups such that the symmetric elements in the group algebras commute, \cite[Remark 3, p. 451]{GSV:98}.
We recall that a non-abelian group $G$ is a Hamiltonian group if every subgroup of $G$ is normal. It is well-known that in this
case $G\cong \mathcal{Q}_8\times E\times O$, \cite[Theorem 1.8.5, p. 63]{PS:02}, where $\mathcal{Q}_8=\langle x, y: x^4=1, x^2=y^2, y^{-1}xy=x^{-1}\rangle$
is the quaternion group of order $8$, $E$ is an elementary abelian $2$-group and $O$ is an abelian group with every element of
odd order. When $O=\{1\}$, $G$ is called a Hamiltonian $2$-group.

When one works with linear extensions of arbitrary involutions of the base group of the group algebra, $\ast: G \rightarrow G$,
one finds a larger class of groups such that the symmetric elements also commute. We recall that a group $G$ is said to be an
LC-group (a group with ``limited commutativity'' property) if it is non-abelian and for any pair of elements $g, h\in G$, we have
that $gh=hg$ if and only if at least one element of $\{g, h, gh\}$ lies in $\zeta(G)$, where $\zeta(G)$ denotes the center of $G$.
This family of groups was introduced by Goodaire. From \cite[Proposition III.3.6, p. 98]{GJP:96}, a group $G$ is an LC-group with
a unique non-trivial commutator $s$ (which must have order $2$ and be central) if and only if $G/\zeta(G)\cong C_2\times C_2$, where
$C_2$ is the cyclic group of order $2$. If $G$ is endowed with an involution $\ast$, then we say that $G$ is a special LC-group, or
SLC-group, if it is an LC-group, it has a unique non-trivial commutator $s$ and on such a group, the map $\ast$ is defined by

\begin{equation}\label{eq1}
g^* = \begin{cases}
         g, & \text{if $g$ is central;} \\
         sg, & \text{otherwise.}
      \end{cases}
\end{equation}

We refer to this as the \emph{canonical} involution on an SLC-group. For instance, if $\ast$ is the classical
involution, then  from expression \eqref{eq1} all elements have order $1$, $2$, or $4$.
Furthermore, if $g$ is a non-central element, then $g^2=s$ and we obtain that every cyclic subgroup of $G$ is
normal, and thus in this case the SLC-groups are precisely the Hamilto\-nian $2$-groups.\\

When we consider on $\mathbb{F}G$ the {\it oriented group involution}
$\Big(\sum_{g\in G} \alpha_gg\Big)^{\circledast} = \sum_{g\in G} \alpha_g\sigma(g)g^{\ast}$, where $G$ is a group with a non-trivial
homomorphism $\sigma:G\rightarrow \{\pm1\}$ and an involution $\ast$, the kernel of $\sigma$ is a subgroup $N$ in $G$ of index $2$.
It is clear that the involution $\circledast$ coincides on the subalgebra $\mathbb{F}N$ with the group involution $\ast$. Also, we
have that the symmetric elements in $G$, under $\circledast$, are the symmetric elements in $N$ regarding $\ast$. If we denote the
sets of symmetric elements in $G$, under the involutions $\circledast$ and $\ast$, by $N^+$ and $G^+$, respectively, then we can
write $N^+=N\cap G^+$. In recent years, this type of involution has been of interest and some results were obtained in the study
of properties of $\mathbb{F}G^+, \mathbb{F}G^-$ and $\mathcal{U}^+(\mathbb{F}G)$, see \cite{OP,CP12,HC:20}. For instance, the authors
in \cite{HC:20} proved that $\mathbb{F}G$ satisfies the $\circledast$-PI, i.e., $\mathbb{F}G$ satisfies a PI where $x_i^{\circledast}$
for some $i$'s appear,  $\alpha^\circledast\alpha=\alpha\alpha^\circledast$ if and only if the set $\mathbb{F}G^+$ of symmetric elements
in regard to $\circledast$ is commutative. Since $[G:N]=2$, the structure of the group $N$ and the action of $*$ on $N$ are both known,
see \cite[Theorem 2.4, p.730]{JRM:06} and \cite[Theorem 3.1, p. 4395]{HC:20}, then this classification depend on whether $N=ker(\sigma)$
is either abelian or an SLC-group. However, this result does not provide a complete description of $G$ and the action of $\circledast$
on $G$. Therefore this is the principal aim in this kind of research.\\

In this paper, we extend the results obtained by Dooms and Ruiz \cite{DRM:07} to the case of the oriented group involution \eqref{eq0}.
More precisely, we classify under middle hypothesis, the groups with a regular group algebra over an infinite field $\mathbb{F}$ of
$\car(\mathbb{F})\neq 2$ for which the $\circledast$-symmetric units satisfy a GI. Further, we prove that if the $\circledast$-symmetric
units of $\mathbb{F}G$, where $G$ is a locally finite group, satisfies a group identity then $\mathbb{F}G$ satisfies a polynomial identity,
see Theorem \ref{teor1} and Theorem \ref{teor2}. Moreover, in the case when the prime radical $\eta(\mathbb{F}G)$ of $\mathbb{F}G$
is nilpotent we characterize the groups for which the symmetric units $\mathcal{U}^+(\mathbb{F}G)$ do satisfy a group identity,
see Theorem \ref{teor3}.\\

Throughout this paper $\mathbb{F}$ will always denote an infinite field with $\car(\mathbb{F})\neq 2$, $G$ a group, $\ast$ and $\sigma$
an involution and a non-trivial orientation of $G$, respectively. We will denote with $\circledast$ an oriented group involution of
$\mathbb{F}G$ given by expression \eqref{eq0}, which is linear extension of the involution $\ast$ of $G$, twisted by the homomorphism $\sigma$.

%%%%%%%%%%%%%%%%%%%%%%%%%%%%%%%%%%%%%%%%%%%%%%%%%%%%%%%%%%%%%%%%%%%%%%%%%%%
\section{Preliminaries and notations}

Let $R$ be a ring with involution $\star$. Let $\mathcal{U}(R)$ its group of units and $\mathcal{U}^+(R)=\mathcal{U}(R)\cap R^+$ the set of
symmetric units. It is well-known that central idempotents are very important in the study of group identities. Moreover the following fact
is proved.

\begin{lem}[{\cite[Theorem 2]{GSV:98}}]\label{lem1}
Let $R$ be a semiprime ring with involution $\star$ such that $\mathcal{U}^+(R)$ is GI. Then every symmetric idempotent
of $R$ is central.
\end{lem}

For a given prime $p$, an element $x\in G$ will be called a $p$-element if its order is a power of $p$ and it is called $p'$-element
if its order is finite and, not divisible by $p$. Moreover, a torsion subgroup $H$ of $G$ is a $p'$-subgroup if every element $h\in H$
is a $p'$-element. We agree that if $p=0$ every torsion subgroup is a $p'$-subgroup.

An immediate consequence in the setting of group algebras of Lemma \ref{lem1} is the following: Let $\mathbb{F}$ be
a field with $\car(\mathbb{F})\geq 0$ and $G$ a group such that $\mathbb{F}G$ is semiprime. If $\mathcal{U}^+(\mathbb{F}G)$
is GI under the classical involution, then every torsion $p'$-subgroup of $G$ is normal in $G$. Indeed to prove this result is important
the use of the idempotent element

$$
\frac{1}{o(g)}\widehat{g}=\frac{1}{o(g)}(1+g+\cdots+g^{o(g)}),
$$

where $g$ is a $p'$-element, see the proof in \cite[Corollary 2, p. 451]{GSV:98}. Note that this element is symmetric under the classical
involution, but it is not when $\ast$ is a group involution.  This fact, in the former case,  gives important information on cyclic subgroups.
Unfortunately, this  property is lost in the case of a group involution $\ast$.

Semisimple algebras whose units satisfy a GI were widely studied, see for instance \cite[Theorem 1, p. 197]{GP:03},
\cite[Theorem 2.2, p. 743]{DRM:07}, \cite[Theorem 3.1, p. 1732]{BDR:09} and \cite[Lemma 2.1, p. 2803]{GPS:09i}.
The following three lemmas will be needed in the sequel.

\begin{lem}[{\cite[Lemma 2.1]{GPS:09i}}]\label{lem2}
Let $R$ be a finite dimensional semisimple algebra with involution $\star$ over an infinite field $K$, $\car(K)\neq 2$. Suppose
that $\mathcal{U}^+(R)$ is GI. Then $R$ is a direct sum of simple algebras of dimension at most four over their centers
and the symmetric elements $R^+$ are central in $R$, i.e.,
$$R\cong D_1\oplus D_2\oplus\cdots \oplus D_k\oplus M_2(\mathbb{F}_1)\oplus M_2(\mathbb{F}_2)\oplus\cdots\oplus M_2(\mathbb{F}_l)~ \text{ and } ~R^+\subseteq \zeta(R).$$
\end{lem}

\begin{lem}[{\cite[Theorem 2.2]{DRM:07}}]\label{lem3}
Let $R$ be a semisimple $K$-algebra with involution $\star$, where $K$ is an infinite field with $\car(K)\neq 2$. Suppose one of the following conditions holds:
 \begin{enumerate}[\normalfont(i)]
 \item $K$ is uncountable,
 \item $R$ has no simple components that are non-commutative division algebras other than quaternion
                    algebras.
\end{enumerate}
Then $\mathcal{U}^+(R)\in \text{GI}$ if and only if $R^+$ is central in $R$.
\end{lem}

\begin{lem}[{\cite[Lemma 2.3.5]{Lee:10}}]\label{lem4}
Suppose that $R$ is an $K$-algebra with involution $\star$, where $\car(K)\neq 2$. Let $I$ be a $\star$-invariant nil ideal.
If $\mathcal{U}^+(R)$ satisfies the group identity $\omega(x_1,\ldots,x_n)=1$, then so does $\mathcal{U}^+(R/I)$. Conversely,
if $p>0$, $I$ is nil of bounded exponent and $\mathcal{U}^+(R/I)\in \text{GI}$, then $\mathcal{U}^+(R)\in \text{GI}$.
\end{lem}

\begin{obs}
By \cite[Lemma 2.4, p. 891]{GPS:09}, under the assumptions of Lemma \ref{lem3}, it follows that $\mathcal{U}^+(R)\in \text{GI}$ is
equivalent to $R^+$ is Lie $n$-Engel, for some $n$.
\end{obs}

We conclude this section with a result due to Jespers and Ruiz Mar\'in \cite{JRM:06}, where the SLC groups arise naturally.

\begin{lem}[{\cite[Theorem 2.4]{JRM:06}}]\label{lem5}
Let $\mathbb{F}$ be a field with $\car(\mathbb{F})\neq 2$ and let $G$ be a group with an involution $\ast$ extended
$\mathbb{F}$-linearly to $\mathbb{F}G$. Then $\mathbb{F}G^+$ is commutative if and only if $G$ is abelian or an SLC group.
In this case, $\mathbb{F}G^+=\zeta(\mathbb{F}G)$.
\end{lem}

%%%%%%%%%%%%%%%%%%%%%%%%%%%%%%%%%%%%%%%%%%%%%%%%%%%%%%%%%%%%%%%%%%%%%%%%%%%
\section{Results}
\subsection{Regular group algebras}

We need the following results about group algebras endowed with an oriented group involution. We recall that for a fixed orientation $\sigma$ of $G$, we denote with $N=\ker(\sigma)$.

\begin{lem}[{\cite[Lemma 1.1]{OP}}]\label{lem7}
Let $R$ be a commutative ring with unity of characteristic different from $2$  and let $G$ be a group with a non-trivial
orientation $\sigma$ and an involution $\ast$. Suppose that $RG^+$ is commutative under oriented group involution and let
$g\in (G\setminus N)\setminus G^+$, $h\in G$. Then one of the following holds:
 \begin{enumerate}[\normalfont(i)]
 \item $gh=hg$; or
 \item $\car(R)=4$ and $gh=g^\ast h^\ast=hg^\ast=h^\ast g$.
\end{enumerate}
Furthermore, $gg^\ast=g^\ast g$.
\end{lem}

\begin{lem}[{\cite[Theorem 2.2]{OP}}]\label{lem8}
Let $R$ be a commutative ring with unity of characteristic different from $2$  and let $G$ be a non-abelian group with
a non-trivial orientation $\sigma$ and an involution $\ast$. Then, $RG^+$ is a commutative ring if and only if one of
the following conditions holds:
 \begin{enumerate}[\normalfont(i)]
 \item $N=ker(\sigma)$ is an abelian group and $(G\setminus N)\subset G^+$;
 \item $G$ and $N$ have the \emph{LC}-property, and there exists a unique non-trivial commutator $s$ such
                    that the involution $\ast$ is given by
\begin{equation*}
g^* = \begin{cases}
         g, & \text{if $g\in N\cap\zeta(G$) or $g\in (G\setminus N)\setminus\zeta(G)$;} \\
         sg, & \text{otherwise.}
      \end{cases}
\end{equation*}
 \item[\emph{(iii)}] $\car(R)=4$, $G$ has the \emph{LC}-property, and there exists a unique non-trivial commutator $s$
                    such that the involution $\ast$ is the canonical involution.
\end{enumerate}
\end{lem}

Recall that a ring $R$ with identity is said to be (von Neumann) regular if for any $x\in R$ there exists an $y\in R$ such that $xyx=x$.
Villamayor \cite[Theorem 3.1.5, p. 69]{Pas:77} showed that the group algebra $\mathbb{F}G$ is regular if and only if $G$ is locally finite
and has no elements of order $p$ in case $\mathbb{F}$ has characteristic $p>0$. Note that in this case the set of $p$-elements $P$ is trivial
and thus $\mathbb{F}G$ is semiprime, \cite[Theorem 4.2.13, p. 131]{Pas:77} (in case $\car(\mathbb{F})=0$, we agree that $P=\left\{1\right\}$).

We are now able to classify the groups with a regular group algebra over an infinite field $\mathbb{F}$ of $\car(\mathbb{F})\neq 2$ for which
the $\circledast$-symmetric units satisfy a GI, result which is the oriented version of Dooms and Ruiz, see Theorem \ref{teor01}.
Note that in this context the third condition in Lemma \ref{lem8} will not be considered.

\begin{teor}\label{teor1}
Let $\mathbb{F}$ be an infinite field with $\car(\mathbb{F})\neq 2$ and let $G$ be a non-abelian group such that
$\mathbb{F}G$ is regular. Let $\sigma:G\rightarrow \{\pm 1\}$ be a non-trivial orientation and an involution $\ast$ on $G$.
Suppose one of the following conditions holds:

 \begin{enumerate}[\normalfont(i)]
 \item $\mathbb{F}$ is uncountable,
 \item all finite non-abelian subgroups of $G$ which are $\ast$-invariant have no simple components in their group algebra over $\mathbb{F}$
       that are non-commutative division algebras other than quaternion algebras.
\end{enumerate}
Then $\mathcal{U}^+(\mathbb{F}G)\in \text{GI}$ if and only if one of the following conditions holds:
 \begin{enumerate}[\normalfont(1)]
 \item $N=ker(\sigma)$ is an abelian group and $(G\setminus N)\subset G^+$;
 \item $G$ and $N$ have the \emph{LC}-property, and there exists a unique non-trivial commutator $s$ such that
       the involution $\ast$ is given by
\begin{equation}\label{eq2}
g^* = \begin{cases}
         g, & \text{if $g\in N\cap\zeta(G$) or $g\in (G\setminus N)\setminus\zeta(G)$;} \\
         sg, & \text{otherwise.}
      \end{cases}
\end{equation}
\end{enumerate}
Consequently, $\mathcal{U}^+(\mathbb{F}G)\in \text{GI}$ if and only if $\mathcal{U}^+(\mathbb{F}G)$ is an abelian
group.
\end{teor}
\begin{proof}
Assume that $\mathcal{U}^+(\mathbb{F}G)\in \text{GI}$ and let $N=ker(\sigma)$. Then $\mathcal{U}^+(\mathbb{F}N)\in \text{GI}$. Hence, by Theorem \ref{teor01} and Lemma \ref{lem5}, we have two possibilities for $N$; either

 \begin{enumerate}[\normalfont(A)]
\item $N$ is an abelian group; or

\item $N$ has the \emph{LC}-property, and there exists a unique non-trivial commutator $s$ such that the involution
           $\ast$ in $N$, is the canonical involution.
\end{enumerate}
Now, let $g, h\in G$ such that $gh\neq hg$. Consider the collection $\{H_i\}_{i\geq 1}$ (possibly infinite) of all finite subgroups of $G$ which
are $\ast$-invariant and contain $g$ and $h$. It is clear that $G=\bigcup_iH_i$ and that $\mathcal{U}^+(\mathbb{F}H_i)\in \text{GI}$. Since $\mathbb{F}G$ is regular, we have that all $\mathbb{F}H_i$ are semisimple.

Set $\sigma_i=\sigma|_{H_i}$ and let $N_i=ker(\sigma_i)$. By Lemma \ref{lem3}, $\mathbb{F}H_i^+$ is central in
$\mathbb{F}H_i$ for all $i$, and applying the Lemma \ref{lem8}, one of the following conditions holds:

 \begin{enumerate}[\normalfont(a)]
\item $N_i$ is an abelian group and $(H_i\setminus N_i)\subset H_i^+$; or
\item $H_i$ and $N_i$ have the \emph{LC}-property, and there exists a unique non-trivial commutator $s$ such that
           the involution $\ast$ is as given in the expression \eqref{eq2}.
\end{enumerate}

It is easy to see that, $N=ker(\sigma)=\bigcup_i N_i=\bigcup_i ker(\sigma_i)$ and $\bigcup_i(H_i\setminus N_i)=G\setminus N$. Suppose that (A) is true. Then (a) holds for all $i$ and thus (1) follows.

Assume that (B) is true. If there exists $j$ such that (a) holds, then since $N_j$ is abelian and  $g, h\in H_j$, at least one of them belong to $H_j\setminus N_j$. Without loss of generality, if $g\in H_j\setminus N_j$ follows that $g$ is symmetric, and as $\mathbb{F}H_j^+$ is central in
$\mathbb{F}H_j$, we obtain $gh=hg$, which is a contradiction.  So (b) holds for all $i$.

If $s$ is not a unique non-trivial commutator of $G$, then there exist $x, y\in G$ such that $(x, y)\neq s$. We know that
$x, y\in H_i$, for some $i$, for instance, $H_i=\langle x,y,g,h,x^\ast,y^\ast,g^\ast,h^\ast\rangle$. Therefore
$(x,y)=(g,h)=s$, a contradiction.

\emph{\textbf{Claim:}} For all $i$, $\zeta(N_i)=N_i\cap \zeta(H_i)$.
In fact, let $x\in \zeta(N_i)\setminus \zeta(H_i)$. Then, there exists $y\in H_i\setminus N_i$ such that $xy\neq yx$ and by the behavior of $\ast$ on $H_i$ given by \eqref{eq2}, $y^\ast=y$. Since $\mathbb{F}N_i^+$ is commutative, we have that $x^\ast=x$. Thus $xy\in (H_i\setminus N_i)\setminus \zeta(H_i)$ and by (b) $xy\in H_i^+$. Therefore $xy=(xy)^\ast=y^\ast x^\ast=yx$,
a contradiction. Hence $\zeta(N_i)=N_i\cap \zeta(H_i)$.

% we have showed above that for all $i$, $\mathbb{F}H_i^+$ is a commutative ring. Let $h\in H_i\setminus N_i$. If $h$ is not central,  then there exists $g\in H_i$ such that $hg\neq gh$ and by Lemma \ref{lem7}, it follows that $h\in H_i^+$.
%If $h$ is central, take $x\in N_i\setminus \zeta(N_i)$. Then $xh\in (H_i\setminus N_i)\setminus \zeta(H_i)$ and again,
%as above $xh\in H_i^+$. Therefore, $xh=(xh)^\ast=h^\ast x^\ast=sh^\ast x$, i.e., $h^\ast=sh$.

Thus for  all $i$, we obtain that $H_i^+=\zeta(N_i)\cup \left[(H_i\setminus N_i)\setminus \zeta(H_i)\right]$. Now, it is clear that $\zeta(G)\subseteq\bigcup_i\zeta(H_i)$. When $x\in\bigcup_i\zeta(H_i)$, then $x\in \zeta(H_j)$ for
some $j$. By the construction of $H_j$, $g\in H_j$ and hence $xg=gx$. Since $g$ is an arbitrary element of $G\setminus \zeta(G)$, we conclude that $x\in \zeta(G)$ and $\zeta(G)=\bigcup_i\zeta(H_i)$. Therefore $\ast:G\rightarrow G$ is given as in the statement, and (2) holds.

The converse is clear, because conditions $(1)$ and $(2)$ by Lemma \ref{lem8} imply that $\mathbb{F}G^+$ is commutative and hence,
$(u,v)=1$ is a GI for $\mathcal{U}^+(\mathbb{F}G)$.

The last assertion is now clear.
\end{proof}

We can find group algebras that fulfil the condition (2) in Theorem \ref{teor1}, see Remark in \cite[p. 746]{DRM:07} and the references quoted therein. \\

%The next lemma which we need is proved in \cite[Lemma 2.6]{GPS:09}.
%
%\begin{lem}\label{lem9}
%Let $\mathbb{F}$ be a field with $\car(\mathbb{F})=p>2$, $G$ a finite group and $\mathcal{J}$ the Jacobson radical of
%$\mathbb{F}G$. Suppose that $\mathbb{F}G/\mathcal{J}$ is isomorphic to a direct sum of simple algebras of dimension at
%most four over their centers. Then the set $P$ of $p$-elements of $G$ is a subgroup.
%\end{lem}
\subsection{Non-regular group algebras}
Dooms and Ruiz, in \cite[Lemma 3.3, p. 747]{DRM:07}, assuming that $\mathbb{F}$ is an infinite field with $\car(\mathbb{F})\neq 2$ and $G$ a locally finite
group such that $\mathcal{U}^+(\mathbb{F}G)\in \text{GI}$ under a group involution $\ast$, demonstrated the set of $p$-elements of $G$ is a normal subgroup
of $G$. They obtained a similar result to Theorem \ref{teor1} for non-regular group algebras, see \cite[Theorem 3.4, p. 748]{DRM:07}.

To handle group algebras which are not necessarily regular, we need the following two lemmas which are the natural extensions of known results. As usual,
for a normal subgroup $H$ of $G$ we denote by $\Delta(G,H)$ the kernel of the map $\mathbb{F}G\stackrel{\Psi}{\rightarrow} \mathbb{F}(G/H)$ defined by

$$
\sum_{g\in G}\alpha_gg\mapsto \sum_{g\in G}\alpha_ggH
$$

and $\Delta(G,G)=\Delta(G)$ is the augmentation ideal.

\begin{lem}\label{lem10}
Let $G$ be a locally finite group and $\car(\mathbb{F})=p\neq 2$. If $\mathcal{U}^+(\mathbb{F}G)\in \text{GI}$, then the set $P$ of $p$-elements of $G$ is a
subgroup.
\end{lem}
\begin{proof}
Let $g,h \in P$ and let $H=\langle g,h,g^\ast,h^\ast\rangle$. Since $G$ is locally finite, then $H$ is finite. Moreover, $H$ is $\ast$-invariant and
$H\subset N=ker(\sigma)$ (every element $x\in H$ has odd order). Since $\mathbb{F}H$ is a finite dimensional algebra the Jacobson radical $\mathcal{J}$
is nilpotent. Let $R=\mathbb{F}H/\mathcal{J}$. Then $R$ is semisimple  and, by Lemma \ref{lem4}, $\mathcal{U}^+(R)$ satisfies a group identity. Hence by
Lemma \ref{lem2}, $R$ is a direct sum of simple algebras of dimension at most four over their centers. Finally, by \cite[Lemma 2.6, p. 892]{GPS:09} we get
that $P$ is a subgroup.
\end{proof}

\begin{lem}\label{lem11}\label{lem13}
Let $\mathbb{F}$ be a field with $\car(\mathbb{F})=p>2$ and $G$ a group such that $\mathcal{U}^+(\mathbb{F}G)$ satisfies a GI $\omega(x_1,...,x_n)=1$,
under an oriented group involution $\circledast$. If $H$ is a normal $\ast$-invariant $p$-subgroup of $G$, and either $H$ is finite or $G$ is locally finite,
then $\mathcal{U}^+(\mathbb{F}(G/H))$ satisfies $\omega(x_1,...,x_n)=1$.
\end{lem}
\begin{proof}
Since $H$ is a $p$-subgroup, then $H\subset N$ and hence $H$ is $\circledast$-invariant.
\\*
If $H$ is finite, then by \cite[Lemma 1.1.1, p. 1]{Lee:10}, $\Delta(G,H)$ is nilpotent and the statement follows from Lemma \ref{lem4}.

Now assume that $G$ is locally finite. Let $\overline{G}=G/H$ and take $\overline{\alpha_1},\ldots,\overline{\alpha_n}\in \mathcal{U}^+(\mathbb{F}\overline{G})$.
Since the map $\mathbb{F}G\rightarrow \mathbb{F}\overline{G}$ is an epimorphism we have that $\mathbb{F}\overline{G}^+$ is the image of $\mathbb{F}G^+$, thus we
may lift these elements up to $\alpha_1, \alpha_2,\ldots,\alpha_n\in \mathbb{F}G^+$ and similarly for their inverses. Let $G_1$ the subgroup of $G$ generated by
the supports of all of these elements. As $G$ is locally finite, $G_1$ is finite. Taking $H_1=G_1\cap H$, we have by the finite case that
$\mathcal{U}^+(\mathbb{F}(G_1/H_1))$ satisfies $\omega(x_1,...,x_n)=1$. Replacing $G$ with $G_1$ and $H$ with $H_1$, we get the statement.
\end{proof}

Consider the group algebra $\mathbb{F}G$, where $G$ is a locally finite group with an oriented group involution $\circledast$, $P$ the set of $p$-elements of $G$
and $\mathbb{F}$ an infinite field with $\car(\mathbb{F})=p\neq 2$. Suppose that $\mathcal{U}^+(\mathbb{F}G)\in \text{GI}$ and let $\overline{G}=G/P$. If
$\mathbb{F}\overline{G}$ satisfies conditions (i) and (ii) of Theorem \ref{teor1}, then by Lemma \ref{lem10}, we have that $P$ is a normal subgroup of $G$ and by
Lemma \ref{lem13} $\mathcal{U}^+(\mathbb{F}\overline{G})$ is GI. Since $\mathbb{F}\overline{G}$ is regular, it follows that either $\overline{G}$ is abelian, or
$\overline{N}$ and $\overline{G}$ satisfy one of the conclusions of Theorem \ref{teor1} and the involution $\overline{\ast}:\overline{G}\rightarrow \overline{G}$
is as given in the expression \eqref{eq2}. Moreover, $\mathbb{F}\overline{G}^+$ is central in $\mathbb{F}\overline{G}$ and thus $\mathbb{F}\overline{G}$ is PI.
Since $\mathbb{F}\overline{G}\cong \mathbb{F}G/\Delta(G,P)$ and $\Delta(G,P)$ is a nil subring of $\mathbb{F}G$ $\circledast$-invariant, by
\cite[Remark 2, p. 450]{GSV:98} we have that $\Delta(G,P)$ is PI and as being PI is closed under ideal extensions, we get that $\mathbb{F}G$ is also PI. Therefore,
we obtain the next result for non-regular group algebras.

\begin{teor}\label{teor2}
Let $g\mapsto g^\ast$ be an involution on a locally finite group $G$, $\sigma:G\rightarrow \{\pm 1\}$ a non-trivial orientation with $N=ker(\sigma)$, $P$ the set
of $p$-elements of $G$ and $\mathbb{F}$ an infinite field with $\car(\mathbb{F})=p\neq 2$. Suppose that $\mathcal{U}^+(\mathbb{F}G)\in \text{GI}$ and that one of
the following conditions holds:
 \begin{enumerate}[\normalfont(i)]
 \item[\emph{(i)}] $\mathbb{F}$ is uncountable,
 \item[\emph{(ii)}] all finite non-abelian subgroups of $G/P$ which are $\ast$-invariant have no simple components in their group algebra over $\mathbb{F}$ that are
                    non-commutative division algebras other than quaternion algebras,
\end{enumerate}
then we have that
 \begin{enumerate}[\normalfont(1)]
 \item $\overline{G}=G/P$ is abelian, or

 \item $\overline{N}=N/P=ker(\overline{\sigma})$ is abelian and $(\overline{G}\setminus \overline{N})\subset
       \overline{G^+}$, or

 \item $\overline{G}$ and $\overline{N}$ have the \emph{LC}-property and there exists a unique non-trivial commutator
       $\overline{s}$ such that the involution $\overline{\ast}$ in $\overline{G}$ is given by

\begin{equation}
\overline{g^*} = \begin{cases}
                  \overline{g}, & \text{if $\overline{g}\in \overline{N}\cap\zeta(\overline{G})$ or $\overline{g}\in (\overline{G}\setminus \overline{N})\setminus\zeta(\overline{G})$;} \\
                  \overline{sg}, & \text{otherwise.}
                 \end{cases}
\end{equation}
Moreover, $\mathbb{F}G\in \text{PI}$.
\end{enumerate}
\end{teor}

To obtain sufficient conditions for locally finite groups $G$ such that $\mathcal{U}^+(\mathbb{F}G)\in \text{GI}$,
we need the following lemma:

\begin{lem}\label{lem12}
Let $\mathbb{F}$ be a field with $\car(\mathbb{F})=p\neq 2$. Let $G$ be a locally finite group with an involution $\ast$ and a non-trivial orientation $\sigma$.
If $P$ is a subgroup of bounded exponent, and either $G/P$ is abelian, or $G/P$ and $N/P$ are as in Theorem \ref{teor1}, then $\mathcal{U}^+(\mathbb{F}G)\in \text{GI}$.
\end{lem}
\begin{proof}
Suppose $P$ is a subgroup of bounded exponent and that $N/P$, $G/P$, $\ast$ and $\sigma$ are as in the statement. Then by Lemma \ref{lem8},
$\mathcal{U}^+(\mathbb{F}\overline{G})$ is abelian. Hence $(\mathcal{U}^+(\mathbb{F}G),\mathcal{U}^+(\mathbb{F}G))\subset 1+\Delta(G,P)$. Now $\Delta(G,P)$ is nil of
bounded exponent and thus $(\mathcal{U}^+(\mathbb{F}G),\mathcal{U}^+(\mathbb{F}G))^{p^n}=1$ for some $n\geq 0$. Hence $\mathcal{U}^+(\mathbb{F}G)\in \text{GI}$.
\end{proof}

\begin{obs}
Note that, under the assumptions of Lemma \ref{lem12}, in case $\mathbb{F}G$ is PI and $G/P$ is abelian, we obtain that $G'\subseteq P$ is of bounded exponent. Hence
by \cite[Theorem 1.1, p. 657]{Pas:97} even $\mathcal{U}(\mathbb{F}G)$ satisfies a GI. Now, if $G/P$ is an SLC-group one easily deduces that also in this case $G'$	is
of bounded exponent, but not necessarily a $p$-group. Finally, if $\overline{N}=N/P$ is abelian and $(\overline{G}\setminus \overline{N})\subset\overline{G^+}$
we can not assure that $G'$ is neither of bounded exponent nor a $p$-group.
\end{obs}

In the sequel, we characterize the groups for which the symmetric units $\mathcal{U}^+(\mathbb{F}G)$ under $\circledast$ do satisfy a group identity. For $g\in G$, let
$C_G(g)=\left\{h\in G: hg=gh\right\}$ be the centralizer of $g$ in $G$. Set $\Phi(G)=\left\{ g\in G: \left[G:C_G(g)\right]<\infty\right\}$ the finite conjugacy subgroup
of $G$ and $\Phi_p(G)=\langle P\cap \Phi(G)\rangle$.

\begin{teor}\label{teor3}
Let $g\mapsto g^\ast$ be an involution on a locally finite group $G$, $\sigma:G\rightarrow \{\pm 1\}$ a non-trivial
orientation and $\mathbb{F}$ an infinite field with $\car(\mathbb{F})=p\neq 2$. Suppose that the prime radical
$\eta(\mathbb{F}G)$ of $\mathbb{F}G$ is a nilpotent ideal and that one of the following conditions holds:
 \begin{enumerate}[\normalfont(i)]
 \item $\mathbb{F}$ is uncountable,
 \item all finite non-abelian subgroups of $G/P$ which are $\ast$-invariant have no simple components in their group algebra over $\mathbb{F}$ that are non-commutative
       division algebras other than quaternion algebras.
\end{enumerate}
Then $\mathcal{U}^+(\mathbb{F}G)\in \text{GI}$ if and only if $P$ is a finite normal subgroup and $G/P$ is
abelian or $G/P$ and $N/P$ are as in Theorem \ref{teor1}.
\end{teor}
\begin{proof}
Suppose that $\mathcal{U}^+(\mathbb{F}G)\in \text{GI}$, then by Lemma \ref{lem10} we have that $P$ is a normal subgroup. Now, by Theorem \ref{teor2}, either $G/P$ is
abelian or $G/P$ and $N/P$ are as in Theorem \ref{teor1} and hence, $\mathbb{F}G\in \text{PI}$. Thus by \cite[Theorem 5.2.14, p. 189]{Pas:77}, $\left[G: \Phi(G)\right]<\infty$
and $\left|\Phi'(G)\right|<\infty$. Since $\eta(\mathbb{F}G)$ is nilpotent \cite[Theorem 8.1.12, p. 311]{Pas:77} gives that $\Phi_p(G)=P\cap \Phi(G)$ is a finite normal
$p$-subgroup. As $P\Phi(G)/\Phi(G)\cong P/P\cap \Phi(G)$ is finite, then $P$ is finite.

Now, the converse is clear by Lemma \ref{lem12}.
\end{proof}

As a corollary following the arguments in Dooms and Ruiz \cite[Corollary 3.7, p. 749]{DRM:07}, we obtain a characterization of the locally finite groups with semiprime group
algebras such that the set of $\circledast$-symmetric units is GI.

\begin{cor}
Let $g\mapsto g^\ast$ be an involution on a locally finite group $G$, $\sigma:G\rightarrow \{\pm 1\}$ a non-trivial orientation and $\mathbb{F}$ an infinite field with
$\car(\mathbb{F})=p\neq 2$ such that $\mathbb{F}G$ is semiprime. Suppose one of the following conditions holds:

 \begin{enumerate}[\normalfont(i)]
 \item $\mathbb{F}$ is uncountable,
 \item all finite non-abelian subgroups of $G/P$ which are $\ast$-invariant have no simple components in their group algebra over $\mathbb{F}$ that are non-commutative
       division algebras other than quaternion algebras.
\end{enumerate}
Then $\mathcal{U}^+(\mathbb{F}G)\in \text{GI}$ if and only if $\mathcal{U}^+(\mathbb{F}G)$ is an abelian group.
\end{cor}
\begin{proof}
Suppose $\mathcal{U}^+(\mathbb{F}G)\in \text{GI}$, then following the lines of the proof of Theorem \ref{teor3}, $\mathbb{F}G$ is semiprime PI. Hence by
\cite[Theorem 4.2.13, p. 131]{Pas:77} $\Phi_p(G)=\left\{1\right\}$ and by the previous proof, we get that $P=\left\{1\right\}$. Thus $\mathbb{F}G$ is regular, and the
result follows from Theorem \ref{teor1}.
\end{proof}

\section*{Acknowledgements}
This paper was written while the authors visited the Universidad Industrial de Santander and Universidad de Nari\~no, and they thank the members of these institutions
for their warm hospitality. The authors are also grateful to Professors Antonio Giambruno and César Polcino Milies for their helpful suggestions and useful comments
some years ago at the beginning of this work. A. Holgu\'{i}n-Villa was partially supported by Decanatura Facultad de Ciencias at Universidad Industrial de Santander
and J. Castillo was partially supported by Vicerrector\'{i}a de Investigaciones e Interacci\'on Social at Universidad de Nari\~no.

\bibliographystyle{plain}

\newcommand*{\doi}[1]{\href{http://dx.doi.org/#1}{doi: #1}}

\end{document}